\documentclass{amsart}

\newtheorem{theorem}{Theorem}[section]

\newtheorem{thm}[theorem]{Theorem}
\newtheorem{lem}[theorem]{Lemma}
\newtheorem{prop}[theorem]{Proposition}
\newtheorem{cor}[theorem]{Corollary}
\theoremstyle{definition}

\theoremstyle{remark}

\numberwithin{equation}{section}

\input{tcilatex}

\begin{document}
\title[The Mean]{The mean: axiomatics, generalizations, applications}
\author[J. E. Gray]{John E. Gray}
\address{Naval Surface Warfare Center Dahlgren\\
Dahlgren Division\\
Sensor Technology Branch Q31\\
18444 Frontage Road Suite 327\\
Dahlgren, VA 22448-5161}
\email{John.E.Gray@navy.mil}
\author[A. Vogt]{Andrew Vogt}
\address{Department of Mathematics and Statistics\\
Georgetown University\\
Washington DC 20057-1233}
\email{vogta@georgetown.edu}

\begin{abstract}
We present an axiomatic approach to the mean and discuss  generalizations of
the mean, including one due to Kolmogorov based  on the Weak Law of Large
Numbers. We offer examples and  counterexamples, describe conventional and
unconventional uses of  the mean in statistical mechanics, and resolve an
anomaly in  quantum theory concerning apparent simultaneous coexistence of 
means and variances of observables. These issues all arise from the 
familiar definition of the mean.
\end{abstract}

\maketitle















\section{Introduction}

The most important number summarizing a data set is generally thought to be
the mean. Some have questioned its utility, comparing it unfavorably with
the median, the mode, the midrange. Capitalists and communists used to argue
over whether mean income or median income was the truer measure of citizen
well-being. For another example, see Kosko \cite{Kosko}. The mean is not
robust against outliers: it can be strongly influenced by a single
observation. This is both a strength and a weakness. Kosko objected that not
only does a Cauchy random variable not have a well-defined mean but the
average of independent identically distributed Cauchy random variables is
itself a Cauchy variable with the same distribution and thus averaging does
not reduce variability at all. 
Investigators often pursue the quest for a single number or a small set of
numbers that capture the essence of a data set, make multiple data sets
comparable, and provide order to the world of data sets. As data sets get
larger and larger, thanks to the digital explosion, scrutiny of measures
that compress data becomes more important. Candidates, in addition to those
mentioned above, include entropy and various generalized means, but no one
has arrived at measures clearly superior to the mean and its associated
measure, the root mean squared deviation or standard deviation.

In work with sample data the mean is easy to understand, in contrast with
other notions from probability theory - such as independence, conditional
probability, and even probability itself. Some have argued (e.g., de Finetti 
\cite{deF}, Pollard \cite{Pol}, and Whittle \cite{Whit}) that the mean is
the fundamental notion in probability theory and should occupy the central
place in all treatments of probability.

In this note we review some properties of the mean, consider some
generalizations for cases when the ordinary mean does not exist, and
investigate the significance of the mean in state space theory and quantum
mechanics.

We begin by axiomatizing the notion of sample mean. Along with familiar
axioms for symmetry, homogenity, and translation invariance, we introduce a 
\emph{condensation} axiom that describes the result of replacing arbitrary
values by their sample mean. We then use the Strong Law of Large Numbers to
arrive at the familiar mathematical notion of mean, $E(X)$. 
Thereafter we consider generalizations of the mean. These are not needed for
bounded or semi-bounded random variables, but really only for variables that
have heavy-tailed distributions on both right and left, with tails of
similar size. We consider what happens when a random variable is restricted
to an interval $[c - M, c + M]$ and $M$ is allowed to tend to infinity. We
state a theorem (Theorem 3.1) describing the different kinds of behavior
possible and provide examples of each. One generalization, which is due to
Kolmogorov, is what we have chosen to call the \emph{weak mean}, $E_w(X)$,
and corresponds precisely to validity of the Weak Law of Large Numbers. Yet
another generalization, the \emph{doubly weak mean}, $E_{ww}(X)$, applies to
the Cauchy distribution. We also discuss multipliers that can be applied to
a variable $X$ to finitize the mean in the spirit of Feynman and note the
dangers of such finitizations. Nonetheless, we recognize that attempts to
scrutinize the notion of mean in connection with the Cauchy distribution and
other long-tailed distributions are timely.

Turning to applications, we point out that the mean is a natural tool in
state space theory for the transition from deterministic models to
statistical models. We discuss entropy and observe that although it is
regarded as a mean it is very different from means arising from ordinary
observables. We recall Jaynes' Maximum Entropy Principle, which seeks to
maximize entropy subject to given values of conventional means.

Lastly, we discuss the role the mean plays in quantum theory, and provide a
precise answer to the question of when the mean and variance exist for a
particular quantum state and a particular quantum observable.

The conclusion, implicit in this discussion, is that the mean is the
paramount measure, of great and wide utility, instructive even when it falls
short. There is little prospect of it losing its longtime preeminence.

\section{Axiomatics for the Sample Mean and the Strong Law of Large Numbers}

Prior to introducing probability measures, let us consider potential axioms
for the mean of a finite set. In this setting, with $\mathcal{R} =
(-\infty,\infty)$, the mean can be thought of as a family of functions $%
\{f_{n}\}$ for $n\geq 1$ with $f_{n}:\mathcal{R}^{n}\rightarrow \mathcal{R}$%
. Its properties include the following: \vspace{0.2in}

M-1) (Homogeneity) $f_n (\lambda x_1, \ldots, \lambda x_n) = \lambda
f_n(x_1, \ldots, x_n)$ for all \newline
$(x_1, \ldots, x_n) \in \mathcal{R}^n$ and all $\lambda \in \mathcal{R}$; 
\vspace{.2in}

M-2) (Symmetry) $f_n(x_1, \ldots, x_n) = f_n(x_{\sigma(1)}, \ldots,
x_{\sigma(n)})$ for all permutations $\sigma$ of the set $\{1, 2, ..., n\}$; 
\vspace{.2in}

M-3) (Translation Invariance) $f_n(x_1 + c, \ldots, x_n + c) = f_n(x_1,
\ldots, x_n) + c$ for all $(x_1, \ldots, x_n) \in \mathcal{R}^n$ and all c $%
\in \mathcal{R}$. \vspace{.2in}

Other properties are the following: \vspace{.2in}

(Positive Homogeneity) $f_n (\lambda x_1, \ldots, \lambda x_n) = \lambda
f_n(x_1, \ldots, x_n)$ for all \newline
$(x_1, \ldots, x_n) \in \mathcal{R}^n$ and all $\lambda > 0$; \vspace{.2in}

(Nonnegativity) If for some $(x_1, \ldots, x_n)\mbox{ and } (y_1, \ldots,
y_n) \in \mathcal{R}^n$ \newline
$x_1 \leq y_1, \ldots, x_n \leq y_n$, then $f_n(x_1, \ldots, x_n) \leq
f_n(y_1, \ldots, y_n)$; \vspace{.2in}

(Positivity) If for some $(x_1, \ldots, x_n)\mbox{ and } (y_1, \ldots, y_n)
\in \mathcal{R}^n$ \newline
$x_1 \leq y_1, \ldots, x_n \leq y_n$, and $x_i < y_i$ for some i, then $%
f_n(x_1, \ldots, x_n) < f_n(y_1, \ldots, y_n)$. \vspace{.2in}

(Strict Positivity) If for some $(x_1, \ldots, x_n)\mbox{ and } (y_1,
\ldots, y_n) \in \mathcal{R}^n$ \newline
$x_i < y_i$ for all i = 1, ..., n, then $f_n(x_1, \ldots, x_n) < f_n(y_1,
\ldots, y_n)$. \vspace{.2in}

(Additivity) $f_n(x_1 + y_1, \ldots, x_n + y_n) = f_n(x_1, \ldots, x_n) +
f_n(y_1, \ldots, y_n)$ for all $(x_1, \ldots, x_n)\mbox{ and
} (y_1, \ldots, y_n) \in \mathcal{R}^n$; \vspace{.2in}

The above axioms seem reasonable except for additivity. The measure should
be independent of units, thus homogeneous, and independent of the choice of
zero point, and a function of the set rather than the ordered set. In
addition, to capture characteristics of the data, ordering properties -
nonnegativity and perhaps positivity - are not unreasonable. However,
additivity asserts a relationship between the ordering of two data sets that
survives reordering of one set, and this seems much too restrictive.

Consider a rival measure to the mean, namely, the median. The median of a
finite data set $\{x_1, \ldots, x_n\}$ is defined as the midmost of the
numbers when they are arranged in increasing order if $n$ is odd, and half
the sum of the two midmost numbers in such an arrangement if $n$ is even.

The median satisfies homogeneity, symmetry, translation invariance, and
nonnegativity. Furthermore, any fixed convex combination of the mean and the
median other than the median itself satisfies homogeneity, symmetry,
translation invariance, nonnegativity, positivity, and strict positivity.
Indeed, not only the median, but the maximum and the minimum of $\{x_1,
\ldots, x_n\}$ (and other rank functions and convex combinations) satisfy
positive homogeneity, symmetry, translation invariance, and nonnegativity.%
\vspace{.2in}

\begin{prop}
Let $f_n: \mathcal{R}^n \rightarrow \mathcal{R}$ be a function satisfying
homogeneity, symmetry, and translation invariance. \vspace{.2in}

1) If n = 1, then $f_1(x) = x $ for all $x \in \mathcal{R}$. \vspace{.2in}

2) If n = 2, then $f_2(x_1,x_2) = \frac{x_1 + x_2}{2}$ for all $(x_1,x_2)
\in \mathcal{R}^2$.
\end{prop}

\begin{proof}
Homogeneity implies that $f_1(0) = f_2(0,0) = 0$. Translation invariance
then indicates that $f_1(x) = f_1(0 + x) = f_1(0) + x = x$. When n = 2,

\begin{eqnarray*}
f_2(a, b) & = & f_2( -\frac{b - a}{2} + \frac{a + b}{2}, \frac{b - a}{2} + 
\frac{a + b}{2}) \\
& = & f_2(-\frac{b - a}{2}, \frac{b - a}{2}) + \frac{a + b}{2} \\
& = & (\frac{b - a}{2})f_2(-1, 1) + \frac{a + b}{2}.
\end{eqnarray*}

\noindent However, by homogeneity and symmetry $f_2(-1, 1) = -f_2(1,-1)$ 
\newline
$= -f_2(-1, 1)$, and $f_2(-1,1) = 0$.
\end{proof}

When n = 1 or 2, the median and the mean coincide. However, it is obvious
that they do not coincide in general when n is 3 or larger. Without the
requirement of additivity it is natural to inquire whether there is another
suitable property that will distinguish between the median and the mean. One
property that we consider and reject is that $f_n$ shall have continuous
partial derivatives.

\begin{prop}
Let $f_n: \mathcal{R}^n \rightarrow \mathcal{R}$ be a function satisfying
homogeneity, symmetry, and translation invariance that has partial
derivatives at each point with the partial derivatives continuous at $%
(0,\ldots,0) \in \mathcal{R}^n$. Then

\begin{eqnarray*}
f_n(x_1,\ldots, x_n) = \frac{x_1 + \ldots + x_n}{n}
\end{eqnarray*}

\noindent for all $(x_1,\ldots,x_n) \in \mathcal{R}^n$.
\end{prop}

\begin{proof}
If we differentiate the equation $f_n (\lambda x_1, \ldots, \lambda x_n) =
\lambda f_n(x_1, \ldots, x_n)$ with respect to $x_i$, we obtain:

\begin{eqnarray*}
\lambda \frac{\partial f_n}{\partial x_i}(\lambda x_1, \ldots, \lambda x_n)
& = & \lambda \frac{\partial f_n}{\partial x_i}(x_1, \ldots, x_n).
\end{eqnarray*}

\noindent Cancelling $\lambda$ from each side and taking a limit as $\lambda$
approaches 0, we obtain:

\begin{eqnarray*}
\frac{\partial f_n}{\partial x_i}(0, \ldots, 0) & = & \frac{\partial f_n}{%
\partial x_i}(x_1, \ldots, x_n).
\end{eqnarray*}

\noindent Thus all partial derivatives are constant. Since $f_n(0,\ldots,0)
= 0$ by homogeneity, $f_n$ has the form:

\begin{eqnarray*}
f_n(x_1,\ldots, x_n) & = & a_1x_1 + \ldots + a_nx_n.
\end{eqnarray*}

\noindent Symmetry now dictates that $a_1 = \ldots = a_n$ and the fact that $%
f_n(1, \ldots, 1) = f_n(0,\ldots,0) + 1 = 0 + 1 = 1$ accordingly implies
that each $a_i = \frac{1}{n}$.
\end{proof}

The continuous differentiability assumption seems to be aimed primarily at
elimination of the median. So we reject it. Instead we offer as an axiom a
different property characteristic of the mean. \vspace{.2in}

M-4) (Condensation) For $n > m$,

\begin{eqnarray*}
f_n(x_1, \ldots, x_n) & = & f_n( f_m(x_1, \ldots, x_m), \ldots, f_m(x_1,
\ldots, x_m), x_{m+1}, \dots, x_n)
\end{eqnarray*}

\noindent for all $(x_1, \ldots, x_n) \in \mathcal{R}^n$. \vspace{.2in}

This property asserts that if a subset of data is replaced by its ``mean",
the grand ``mean" is not changed. This is the first property that proposes a
definite relationship between means of sets of different sizes. In view of
the symmetry axiom (M-3), the statement does not really restrict the order
of the subset, and as we shall see shortly the statement is only really
needed in special cases. 
%
%
\vspace{.2in}

\begin{prop}
Let $f_n: \mathcal{R}^n \rightarrow \mathcal{R}$ be a function satisfying
homogeneity, symmetry, translation invariance, and condensation, that is,
M-1, M-2, M-3, and M-4. Then

\begin{eqnarray*}
f_n(x_1,\ldots, x_n) = \frac{x_1 + \ldots + x_n}{n}
\end{eqnarray*}

\noindent for all $(x_1,\ldots,x_n) \in \mathcal{R}^n$.
\end{prop}

\begin{proof}
In view of Proposition 2.1 we need only perform an inductive step showing
that the mean formula holds for $n \geq 3$ when it holds for $n - 1$.
Consider

\begin{eqnarray*}
\lefteqn{f_n (x_1, \ldots, x_n)} \\
& = & f_n (\frac{1}{n - 1}(x_1 + \ldots x_{n-1}), \ldots,\frac{1}{n - 1}(x_1
+ \ldots x_{n-1}), x_n) \\
& = & f_n(0,\ldots, 0, x_n - \frac{1}{n - 1}(x_1 + \ldots x_{n-1})) + \frac{1%
}{n - 1}(x_1 + \ldots x_{n-1}) \\
& = & (x_n - \frac{1}{n - 1}(x_1 + \ldots x_{n-1}))f_n(0, \ldots,0, 1) + 
\frac{1}{n - 1}(x_1 + \ldots x_{n-1}) \\
& = & a_1x_1 + ... + a_nx_n.
\end{eqnarray*}

\noindent This shows that $f_n$ is a linear function of $x_1, \ldots, x_n$.
It now follows from Proposition 2.2 that it is the mean.
\end{proof}

The proof of Proposition 2.3 requires that M-4 holds in the case when $m = n
- 1$. In fact we can get by with the assumption that M-4 holds when $m = 2$.
It is easy to see that in this case M-4 also holds for $m = 2^k$. Now set $n
= 2^k + j$, where $0 \leq j < 2^k$. If $j = 0$, $f_n(x_1,...,x_n) =
f_n(c,...,c) =cf_n(1,..,1)$ where $c = (x_1 + ... + x_n)/n$. If $j > 0$,
then $f_n(x_1,...,x_n) = f_n(c,...,c,x_{m+1},...x_n)$ where $m = 2^k$ and $c
= (x_1 + ... + x_m)/m$. We now replace $x_1,...,x_m$ by $x_1^{\prime
},...,x_{m-j}^{\prime },0, ...,0$ so that $c = (x_1 + ... + x_m)/m =
(x_1^{\prime }+ ... + x_{m-j}^{\prime }+ 0 + ... + 0)/m$. Using the symmetry
and homogeneity axioms, we obtain: $f_n(x_1,...,x_n) = f_n((mc + x_{m+1} +
... + x_n)/m, ...,(mc + x_{m+1} + ... + x_n)/m, 0, ...0) = ((mc + x_{m+1} +
... + x_n)/m)f_n(1,...,1,0,...,0) = (x_1 + ... + x_n)/m)f_n(1.,,,.1,0,...,0)$%
. Thus we have established that $f_n$ is linear in $x_1,...,x_n$. By
Proposition 2.3 $f_n$ is the ordinary mean.

A further note on axiomatics is that the translation invariance axiom can be
replaced by $f_n(1,...,1) = 1$ if we also assume that $f_2(x_1,x_2) = (x_1 +
x_2)/2$.

To pass from the sample mean of a finite set to the usual general notion of
mean, we introduce a real-valued random variable $X$. We suppose that
associated with $X$ is a Borel probability measure $P_{X}$ taking each Borel
subset A of the real numbers to:

\begin{eqnarray*}
P_{X}(A) = \mbox{ the probability that X belongs to the set A}.
\end{eqnarray*}

\noindent The \underline{mean} of X, denoted by E(X) or $\mu_X$, is defined
when $x$ is integrable with respect to $P_{X}$ to be:

\begin{eqnarray*}
E(X) = \int_{\mathcal{R}} x P_X(dx).
\end{eqnarray*}

One direction of the remarkable Strong Law of Large Numbers (see Pollard 
\cite[p. 78 and pp. 37-8]{Pol}) states that if $\{X_n\}$ is a sequence of
independent random variables with common distribution $P_X$ and there exists
a constant $m$ such that

\begin{eqnarray*}
\frac{X_1 + \ldots + X_n}{n} \mbox{ converges almost surely to } m
\end{eqnarray*}

\noindent as $n \rightarrow \infty$, then each $X_n$ has mean $m$. Here
``almost surely'' means outside a set of measure zero in the countably
infinite product space induced by the measure $P_X$ (see \cite[pp. 99-102]%
{Pol}). More briefly, if sample means of independent copies of $X$ settle
down to something, then that something is $E(X)$. This can be regarded as
the motivation for the transition from the sample mean to the mathematical
mean $E(X)$. The general notion of mean is derived from the finitary notion
considered earlier.

The other direction of the Strong Law of Large Numbers asserts that if $E(X)$
exists, then the sample mean of $n$ identical independent copies of $X$
converges almost surely to $E(X)$ as $n$ tend to infinity. For a proof of
both directions of the Strong Law, see \cite[pp. 95-102, p. 105]{Pol}. For
an alternate proof due to N. Etemadi, see \cite[pp. 106-7]{Pol}.

The transition here from finite samples to infinite populations
distinguishes the deductive method from the inductive method. While true
science deals comfortably with induction based on finite samples, the
deductive method of the Greeks (and Isaac Newton) relies on axioms whose
relationship to reality is only approximate and always contingent.


Indeed, in using the Strong Law of Large Number we are admittedly
introducing the full panoply of probability theory. It is possible, as noted
in the Introduction, to represent all of probability theory using the mean
as the primitive notion. Thus $P_X(A)$ can be defined as $E(\chi_A(X))$, the
mean of $\chi_A(X)$, where $\chi_A(X)$ is the random variable that equals 1
when $X$ is in $A$ and $0$ when $X$ is not in $A$. However, since in what
follows we plan to use probability theory in its conventional form (i.e.,
according to the axioms of Kolmogorov \cite{Kol}), we see no reason to
restate measure-theoretic facts in terms of the mean as primitive. Indeed, a
reason not to do so is that the concept of independence, which is also
fundamental in probability, is awkward when expressed exclusively in terms
of means.

\section{Extending the Mean}

When $x$ is not integrable with respect to $P_X$, the notion $E(X)$ above is
inapplicable and we must rely on other notions of mean. Richard Feynman was
famous for his integration tricks, and some of these are recorded in the
book of Mathews and Walker \cite{Mat}, based on lectures Feynman gave at
Cornell. Feynman's tricks partly motivated our investigation.

Perhaps the most obvious generalization is the following:

\begin{eqnarray*}
\mbox{Let } L(c) = \lim_{M \rightarrow \infty} \int_{[c - M, c + M]} \, x
\,P_X(dx)
\end{eqnarray*}

\noindent for a real number $c$.

By the Lebesgue Dominated Convergence Theorem this notion coincides with the
ordinary mean when $x$ is integrable with respect to $P_{X}$. Kolmogorov 
\cite[p. 40]{Kol}, in his great foundational work, noted this option in the
case when $c = 0$ and observed that it does not require integrability of $%
|x| $. Indeed if $X$ is a random variable obeying the Cauchy distribution $%
f(x) = 1/\pi(1 + x^2)$, then $X$ satisfies $L(c) \equiv 0$ for any choice of 
$c$.

We mention two related notions of mean:

\begin{enumerate}
\item[L-1)] $\lim_{M \rightarrow \infty}\int_{[a - M, b + M]} \, xP_{X}(dx) %
\mbox{, and}$

\item[L-2)] $\lim_{\min {\{M,K\}}\rightarrow \infty } \int_{[a - M, b + K]}
\, xP_{X}(dx),$
\end{enumerate}

\noindent where $a \leq b$.

It is easily seen that L-1 coincides with $L((a + b)/2)$ since

\begin{equation*}
[a - M, b + M] = [ \frac{a + b}{2} - (M + \frac{b - a}{2}), \frac{a + b}{2}
+ (M + \frac{b - a}{2})]. 
\end{equation*}

\noindent As for L-2, we have the following result. \vspace{.2in}

\begin{prop}
Let X be a random variable with probability measure $P_{X}$. Then for some $a
$ and $b$ with $a \leq b$

\begin{equation*}
\lim_{\min {\{K,M\}} \rightarrow \infty} \int_{[a-M,b+K]} \, xP_{X}(dx) 
\end{equation*}

\noindent exists if and only if $x$ is integrable with respect to $P_{X}$.
\end{prop}

\begin{proof}
If $x$ is integrable on $\mathcal{R}$, the limit exists and equals the mean
of $x$ by Lebesgue's Dominated Convergence Theorem. Conversely, if the limit
exists, then

\begin{equation*}
0 \leq \int_{(b+K,b+K^{\prime }]} x \, P_{X}(dx) < \epsilon 
\end{equation*}

\noindent for $K < K^{\prime }$, both sufficiently large, and any given $%
\epsilon$. Likewise

\begin{equation*}
-\epsilon < \int_{[a-M^{\prime },a-M)} x \, P_{X}(dx) \leq 0 
\end{equation*}

\noindent for $M < M^{\prime }$, both sufficiently large. Fatou's Lemma or
Levi's Theorem \cite[p. 172]{HS} thus implies that

\begin{equation*}
0 \leq \int_{(b+k,\infty)} x \, P_{X}(dx) \leq \epsilon, \hspace{.2in}
-\epsilon \leq \int_{(-\infty,a-M)} x \, P_{X}(dx) \leq 0 
\end{equation*}

\noindent and thus $x$ is integrable on $[0,\infty)$ as well as $(\infty,0]$
and so is integrable on $\mathcal{R} = (-\infty,\infty)$. \vspace{.2in}

So, when L-2 exists, it coincides with E(X).
\end{proof}

We now return to the study of $L(c)$. We shall allow $-\infty \leq L(c) \leq
\infty$. This gives us a bit more flexibility in characterizing what can
happen.

\begin{lem}
Let X be a random variable with probability measure $P_{X}$, and let $c_1$
and $c_2$ be real numbers with $c_1 < c_2$. Then there are three
possibilities: \vspace{.2in}

\hspace{1 in} i) If $L(c_1)$ exists in $[-\infty, \infty]$, then

\begin{equation*}
L(c_1) \leq \liminf _{M \rightarrow \infty} \int_{[c_2 - M, c_2 + M]} \, x
\,P_X(dx); 
\end{equation*}
\vspace{.2in}

\hspace{1 in} ii) If $L(c_2)$ exists in $[-\infty, \infty]$, then

\begin{equation*}
L(c_2) \geq \limsup _{M \rightarrow \infty} \int_{[c_1 - M, c_1 + M]} \, x
\,P_X(dx); 
\end{equation*}

\vspace{.2in}

\hspace{1in} iii) If $L(c_1)$ and $L(c_2)$ both exist in $[-\infty, \infty]$
, then $L(c_ 1) = L(c_2)$.
\end{lem}

\begin{proof}
Suppose $c_1 < c_2$. Then

\begin{eqnarray*}
\lefteqn{\int_{[c_2 - M,c_2 + M]}\, xP_X(dx) = \int_{[c_1 - M,c_1 + M]}\,
xP_X(dx)} \\
& & + \int_{(c_1 + M,c_2 + M]}\, xP_X(dx) - \int_{[c_1 - M,c_2 - M)}\,
xP_X(dx).
\end{eqnarray*}

\noindent The second and third terms on the right are both non-negative and
accordingly i) and ii) follow. In the case of iii), note that i) and ii)
imply that if both $L(c_1)$ and $L(c_2)$ exist, then $L(c_1) \leq L(c_2)$.

If $L(c_2) - L(c_1) > 0$, then there is a positive constant $K$ (for
example, any positive number $< L(c_2) - L(c_1)$) such that for $M$
sufficiently large:

\begin{eqnarray*}
\lefteqn{ K < \int_{(c_1 + M,c_2 + M]}\, xP_X(dx) - \int_{[c_1 - M,c_2 -
M)}\, xP_X(dx)} \\
& \leq & (c_2 + M)P_X((c_1 + M,c_2 + M]) + (M - c_1)P_X([c_1 - M, c_2 - M))
\\
& \leq & (M + d)(P_X((c_1 + M,c_2 + M] \cup [c_1 - M, c_2 - M))
\end{eqnarray*}

\noindent where $d = \max{\{|c_2|,|c_1|\}}$. Thus

\begin{equation*}
\frac{K}{M + d} < P_X((c_1 + M,c_2 + M] \cup [c_1 - M, c_2 - M)). 
\end{equation*}

\noindent Now replace $M$ by $M_j = M + j(c_2 - c_1)$ for each integer $j
\geq 0$ to get:

\begin{equation*}
\frac{K}{M_j + d} < P_X((c_1 + M_j,c_2 + M_j] \cup [c_1 - M_j, c_2 - M_j)). 
\end{equation*}

\noindent Summing over these inequalities and noting that $c_2 + M_j = c_1 +
M_{j+1}$ and $c_1 - M_j = c_2 - M_{j+1}$, we obtain:

\begin{equation*}
\infty = \sum_{j = 0}^{\infty} \frac{K}{M + d + j(c_2 - c_1)} \leq P_X((c_1
+ M, \infty) \cup (-\infty, c_2 - M)) \leq 1 
\end{equation*}

\noindent for a contradiction. Thus, this case is eliminated. So $L(c_1) =
L(c_2)$.
\end{proof}

\begin{thm}
Let X be a random variable with probability measure $P_{X}$. Then exactly
one of the following possibilities holds:

i) $L(c)$ does not exist in $[-\infty, \infty]$ for any real number $c$;

ii) $L(c)$ exists in $(-\infty, \infty)$ for exactly one real number $c$;

iii) $L(c)$ exists in $[-\infty, \infty]$ for all real numbers $c$ and is
independent of $c$;

iv) there is a number $c_0$ such that $L(c) = \infty$ for $c > c_0$ and $%
L(c) $ does not exist for $c < c_0$; or

v) there is a number $c_0$ such that $L(c) = -\infty$ for $c < c_0$ and $%
L(c) $ does not exist for $c > c_0$.
\end{thm}

\begin{proof}
By Lemma 3.1 it suffices to show what happens when $L(c_2) = L(c_1)$ is
finite. In this case the last two terms in the equation at the beginning of
the proof of Lemma 3.1 each tend to $0$ as $M$ tends to infinity. By a
change of variable, we obtain for the positive number $c = c_2 - c_1$.

\begin{equation*}
\lim_{M \rightarrow \infty} \int_{(M,c+M]} x P_{X}(dx) = \lim_{M \rightarrow
\infty} \int_{[-M-c,-M)} x P_{X}( dx) = 0; 
\end{equation*}

Assume $0 < d < c$ and $M \geq 0$. Then

\begin{equation*}
0 \leq \int_{(M,d+M]} x P_{X}(dx) \leq \int_{(M,c+M]} x P_{X}(dx) 
\end{equation*}

\noindent and

\begin{equation*}
\int_{[-M-c,-M)} x P_{X}( dx) \leq \int_{[-M-d,-M)} x P_{X}( dx) \leq 0 
\end{equation*}

\noindent Thus if ii) holds for $c$, it holds for $d$. On the other hand, if
ii) holds for $c$ it also holds for $nc$ where $n$ is any fixed positive
integer since

\begin{equation*}
\int_{(M,nc+M]} x P_{X}(dx) = \sum_{j = 1}^{n} \int_{((j-1)c+M,jc+M]} x
P_{X}(dx) 
\end{equation*}

\noindent and

\begin{equation*}
\int_{[-M-nc,-M)} x P_{X}( dx) =
\sum_{j=1}^{n}\int_{[-M-(n+1-j)c,-M-(n-j)c)} x P_{X}( dx), 
\end{equation*}

\noindent and if the $M^{\prime }s$ are chosen far enough out so that the
individual integrals are closer to zero than $\epsilon/n$, the sum integral
is within $\epsilon$ of $0$. Finally since any positive real number $d$ is
smaller than $nc$ for some positive integer $n$, all cases are covered.
Accordingly, ii) implies iii).

The argument for i) implies ii) can now be used to show that for any two
real numbers $c_1$ and $c_2$, if either $L(c_1)$ or $L(c_2)$ exists in $%
[-\infty,\infty]$, then the other exists and equals it since the
approximating integrals differ by two integrals on intervals of length $|c_1
- c_2|$ that tend to zero as $M$ tends to infinity. Thus iii) implies iv).
\end{proof}

We give some examples to illustrate that each of the possibilities
enumerated in Theorem 3.1 can occur.

Consider a random variable $X$ whose probability measure is of the form

\begin{equation*}
P_X(A) = \sum_{n = 1}^{\infty} (\frac{1}{2^{2n}} \delta_{2^{2n}}(A) + \frac{1%
}{2^{2n - 1}}\delta_{-2^{2n - 1}}(A)) 
\end{equation*}

\noindent where $\delta_z$ is the (Dirac) probability measure whose value is 
$1$ on any Borel subset $A$ of $\mathcal{R}$ that contains the real number $z
$ and whose value is zero otherwise. The sum of the nonzero values is one,
so this obviously defines a probability measure. However the integral of $x$
over the interval $[c - M, c + M]$ is the difference between the size of the
first set and the size of the second set below:

\begin{eqnarray*}
\, \mbox{ The size of the set } \{n \mbox{ :  } 1 \leq n, 2^{2n} \leq (c +
M) \} & = & \lfloor{\frac{\log{(c+M)}}{2 \log{2}}}\rfloor \\
\mbox{ The size of the set }\{n \mbox{ : } 1 \leq n, 2^{2n-1} \leq (M - c)
\} & = & \lfloor{\frac{\log{(M-c)}}{2\log{2}} + \frac{1}{2}}\rfloor
\end{eqnarray*}

\noindent where $\lfloor \, \,\rfloor$ is the floor function. For fixed $c$
and sufficiently large $M$ the difference of the above quantities can assume
the values $0$ and $-1$ and the integral does not settle down to either one.
This is an instance of Theorem 3.1, part i).

A random variable $X$ can also be defined with probability measure of the
form

\begin{equation*}
P_X(A) = \sum_{n = 1}^{\infty}\frac{1}{2^{n+1}}( \delta_{2^n}(A) +
\delta_{(-2^n)}(A)). 
\end{equation*}

\noindent So $P_X$ is concentrated at the points $\pm 2^n$ and assigns
probability $1/(2^{n+1})$ to such points. For this measure, $L(0)$ equals $0$
by symmetry. However, $L(c)$ does not exist for other choices of $c$. If $c$
is positive, the integral of $x$ over the closed interval $[c - M, c + M]$
reduces to its integral over the open interval $(M - c, M + c]$ and this
integral oscillates between $0$ and $\frac{1}{2}$ for large $M$ depending on
whether $2^n$ is in the interval $(M - c, M + c]$ or not. Similar behavior
occurs when $c < 0$. This example is an instance of Theorem 3.1, part ii).

Now consider a random variable $X$ having a probability density (with
respect to Lebesgue measure on the real line) of the form

\begin{equation*}
f(x) = \left\{ 
\begin{array}{ll}
\frac{1}{1 + Cx^a} \mbox{ if $x \geq 0$} &  \\ 
\frac{1}{1 + D|x|^b} \mbox{ if $x < 0$} & 
\end{array}
\right. 
\end{equation*}

\noindent where $a$ and $b$ are numbers in $(1, 2)$ and $C$ and $D$ are
suitable positive constants that guarantee that the density integrates to $1$%
. Notice that this random variable satisfies iii) of Proposition 3.1. It is
easy to see that $L(c) \equiv \infty$ for all $c$ or $-\infty$ for all $c $
according as $b > a$ or $a > b$. The example illustrates part iii) of
Theorem 3.1 (as does the Cauchy distribution with $L(c) \equiv 0$) .

The probability measure

\begin{equation*}
P_X(A) = \sum_{n = 1}^{\infty} (\frac{2^{n-1}}{3^{n + 1}} \delta_{3^{n}}(A)
+ \frac{2^{n-2}}{3^{n + 1}}\delta_{-3^{n}}(A)) 
\end{equation*}

\noindent illustrates part iv) of Theorem 3.1. If $c \geq 0$, the integral
of $x$ over $[c - M,c + M]$ is given by:

\begin{equation*}
\sum_{\{n \mbox{ :  }1 \leq n, 3^{n} \leq c + M \}} \frac{2^{n - 1}}{3} -
\sum_{\{n \mbox{ :  } 1 \leq n, 3^{n} \leq M - c \}} \frac{2^{n - 2}}{3}, 
\end{equation*}

\noindent and this expression has the value $(2^{n_{0} -1} - 2^{-1})/3$ or $%
(2^{n_{0} -1} + 2^{n_{0} - 2} - 2^{-1})/3$ where $n_{0} \approx (\log{c + M}%
)/\log{3}$ for large $M$. Since $M$ and $n_{0}$ tend to infinity together,
it follows that $L(c) \equiv \infty$ for $c \geq 0$.

On the other hand, if $c = -d$ where $d > 0$, the integral of $x$ over $[c -
M,c+M] = [-M - d, M - d]$ is given by:

\begin{equation*}
\sum_{\{n \mbox{ :  }1 \leq n, 3^{n} \leq M - d\}} \frac{2^{n - 1}}{3} -
\sum_{\{n \mbox{ :  } 1 \leq n, 3^{n}\leq M + d \}} \frac{2^{n - 2}}{3}, 
\end{equation*}

\noindent and this reduces to $(2^{n_{0} -1} - 2^{-1})/3$ or to $(-1)/6$ for
large $M$ depending on whether a positive integer $n_{0}$ lies in the
interval $(\log{(M - d)}/\log{3}, \log{(M + d)}/\log{3}]$ or not. Thus, $L(c)
$ does not exist for $c < 0$.

A final example (also for part iv) of Theorem 3.1 is the case where the
probability measure is given by:

\begin{equation*}
P_X(A) = K \sum_{n = 1}^{\infty} (\frac{2^{n}}{3^{n} + (1/n)} \delta_{3^{n}
+ (1/n)}(A) + \frac{2^{n-1}}{3^{n}}\delta_{-3^{n}}(A)). 
\end{equation*}

\noindent Here $K$ is a suitably chosen positive normalizer, which is easily
seen to be smaller than $1/3$. For $c > 0$, the integral of $x$ over $[c -
M,c + M]$ is

\begin{equation*}
K\sum_{\{n \mbox{ :  }1 \leq n, (3^{n} + 1/n) \leq M + c\}} 2^{n} - K
\sum_{\{n \mbox{ :  } 1 \leq n, 3^{n} \leq M - c \} }2^{n - 1}, 
\end{equation*}

\noindent and this reduces to $K(2^{n_0} - 1)$ for $M$ sufficiently large
where $n_0$ is the largest integer such that $3^{n_0} \leq M - c$. Since $%
n_0 $ and $M$ tend to infinity together, $L(c) \equiv \infty$ for all $c > 0$%
.

When $c = 0$, the integral of $x$ over $[-M,M]$ reduces to $K(2^{n_0} - 1)$
or to $-K$ where $n_0$ is the largest integer such that $3^{n_0} + (1/n_0)
\leq M$ and the first or second reduction occurs according as $M < 3^{n_0 +
1}$ or not. Thus $L(0)$ does not exist. By Theorem 3.1 $L(c)$ does not exist
for $c < 0$.

Other cases arising in Theorem 3.1, such as part v), are obtained by
modifying the examples above, e.g., replacing $X$ by $-X$ or by $X + a$.

\section{Weak Means and Multipliers}

One of the implications of Theorem 3.1 is that if $L(c)$ exists for more
than one choice of $c$ and is finite in some case then it is finite for all $%
c$ and is independent of $c$. The case of the Cauchy distribution shows that
this can happen without the ordinary mean existing. Accordingly for a random
variable $X$, we define the \underline{doubly weak mean} of $X$, denoted by $%
E_{ww}(X)$, to be the common value of $L(c)$ for all $c$ when this common
value exists and is in $(-\infty,\infty)$.

We also introduce an intermediate notion due to Kolmogorov between the
ordinary mean and the doubly weak that motivates our terminology. The 
\underline{weak mean} of $X$, denoted by $E_w(X)$, is defined as follows: $%
E_w(X)$ is the quantity $L(0)$ provided the latter exists in $(-\infty,
\infty)$ and provided $\lim_{n \rightarrow \infty} nP_X(|X| > n) = 0$.

The following proposition is due to Kolmogorov. It indicates that existence
of the weak mean coincides precisely with the existence of a number for
which the Weak Law of Large Numbers holds. \vspace{.2in}

\begin{prop}[Kolmogorov, 1928]
Let $X$ be a random variable. Suppose that $\{ X_1, \cdots, X_n, \cdots \}$
are independent identically distributed copies of $X$ with $P_n$ the n-fold
product distribution. Then there is a real number $m$ such that for each $%
\epsilon > 0$

\begin{equation*}
\lim_{ n \rightarrow \infty} {P_n( | \frac{X_1 + ... + X_n}{n} - m| >
\epsilon)} = 0 
\end{equation*}

if and only if $X$ has weak mean $E_w(X) = m$.
\end{prop}

\begin{proof}
See \cite[p. 65]{Kol}, \cite{Kol2}, and \cite[Theorems XII and XIII]{Kol3}.
\end{proof}

\vspace{.2in} 

\begin{cor}
Let $X$ be a random variable.

i) If $X$ has a mean, then $X$ has a weak mean and $E(X) = E_w(X)$;

and

ii) if $X$ has a weak mean, then $X$ has a doubly weak mean and $E_w(X) =
E_{ww}(X)$.
\end{cor}

\begin{proof}
In case $X$ has a mean, then the identity function $x \mapsto x$ is
integrable with respect to the probability measure $P_X$ on the real line.
In particular the tail integrals

\begin{equation*}
\int_{[n, \infty)} x P_X(dx) \mbox{ and } \int_{(-\infty, -n]} x P_X(dx)
\end{equation*}

\noindent tend to zero as $n$ tends to infinity. Since the absolute values
of these integrals are larger respectively than $nP_X(X \geq n)$ and $nP_X(X
\leq -n)$, it follows that $\lim_{n \rightarrow \infty} nP_X(|X| > n) = 0$.
Likewise by Lebesgue's Dominated Convergence Theorem, $L(0) = E(X)$, 
Suppose $X$ has a weak mean. Then if $c_1 < c_2$ and $\epsilon> 0$ and a
sufficiently large $M$ are given,

\begin{eqnarray*}
0 \leq \int_{(c_1 + M,c_2 + M]}\, xP_X(dx) & \leq & (c_2 + M) P_X((c_1 +
M,c_2 + M]) \\
& \leq & (c_2 - c_1 + c_1 + M) P_X(|X| > c_1 + M) \\
& \leq & (c_2 - c_1 + c_1 + M) P_X(|X| \geq n) \\
& \leq & (\frac{c_2 - c_1}{n} + \frac{c_1 + M}{n})\epsilon
\end{eqnarray*}

\noindent where $n = \lfloor{c_1 + M}\rfloor$. For $M$ sufficiently large,
the right side is as close to $\epsilon$ as we like. Thus

\begin{equation*}
\lim_{M \rightarrow \infty} \int_{(c_1 + M,c_2 + M]}\, xP_X(dx) = 0.
\end{equation*}

\noindent Similarly,

\begin{equation*}
\lim_{M \rightarrow \infty} \int_{[c_1 - M, c_2 - m)}\, xP_X(dx) = 0.
\end{equation*}

\noindent Accordingly from the first equation in the proof of Lemma 3.1, it
follows that when one of $L(c_2)$ or $L(c_1)$ exists and is finite, the
other exists and is equal to it. Since $L(0) = m$, it follows that $L(c)$
exists for all $c$, $L(c) \equiv m$ and $m$ is the doubly weak mean of $X$.
\end{proof}

Kolmogorov in \cite[p. 66]{Kol} gives an example where the Weak Law holds
but the Strong Law does not. Cauchy random variables have $L(c)$ existing
for all $c$, independent of $c$, but violate the Weak Law by not decaying
rapidly enough at infinity. Thus the mean, weak mean, and doubly weak mean
are strictly distinct notions.

We make one more observation on generalizations of the mean, based on using
multipliers to attempt to finitize the mean. These multipliers are a type of
``mollifier." Usually mollifiers are used to aid approximation of the delta
function and to smooth functions, but another use is to regularize behavior
at $\pm \infty$. The idea is to introduce a function $\phi_{\lambda}(x)$
that depends on a parameter $\lambda$ so that $x \mapsto \phi_{\lambda}(x)x$
is integrable with respect to $P_{X}$ for $\lambda \neq \lambda_{0}$ and $%
\phi_{\lambda}(x) \rightarrow 1$ for a. e. x as $\lambda \rightarrow
\lambda_{0}$. In the case of L(c), the multiplier can be taken to be

\begin{equation*}
\phi_{\lambda} (x) = \chi_{[c - 1/\lambda, c + 1/\lambda]}(x) 
\end{equation*}

\noindent where $\chi_A$ is the characteristic function of the set $A$ and $%
\lambda = 1/M$.

Multipliers, and indeed other straight-forward generalizations of the mean
including the weak and doubly weak mean, are useful only when the following
equations hold:

\begin{eqnarray*}
\int_{[0,\infty)} x P_{X}(dx) & = & \infty \\
\int_{(-\infty,0]} x P_{X}(dx) & = & - \infty .
\end{eqnarray*}

\noindent If neither of these equations holds, x is integrable and the mean
is well-defined. If only the first equation holds, the mean is $+\infty $,
and if only the second equation holds, the mean is $-\infty $. If both
equations hold, then there is some room for maneuver. $L(c)$ cannot exist
finitely unless the infinities on each end are of the same order. If for
example $P_{X}$ is given by a density function $f$ with respect to Lebesgue
measure such that $f(x)$ decays as $1/x^{2}$ as $x\rightarrow \infty $ and
decays as $1/|x|^{3/2}$ as $x\rightarrow -\infty $, then $L(c)\equiv -\infty 
$ for all $c$.

Multipliers offer possibilities for extending the notion of the mean. They
can be of use in such activities as renormalization where the aim is to
reinterpret integrals to make them finite. In our case we set:

\begin{equation*}
E_{mult}(X) = \lim_{\lambda \rightarrow \lambda_{0}} {E(\phi_{\lambda}(X)X)} 
\end{equation*}

\noindent provided this limit exists. This method is used to ``evaluate" the
integrals of $\sin{bx}$ and $\sin{x}/x$ on $[0,\infty)$ in \cite[p. 60 and
p. 91]{Mat}.)

However, there are dangers that the following example illustrates.

Define a function $\phi_{\lambda,c}$ for $\lambda > 0$ and $c$ in $\mathcal{R%
}$ by:

\begin{equation*}
\phi_{\lambda,c}(x) = \left\{ 
\begin{array}{ll}
e^{-\lambda x} & \mbox{if $x > 0$} \\ 
e^{\lambda x}(1 + \pi c\lambda x) & \mbox{if $x < 0$.}%
\end{array}%
\right. 
\end{equation*}

\noindent Here $c$ is an arbitrary constant. Evidently $\phi_{\lambda,c}$ is
a well-behaved function, integrable and dying off at $\pm \infty$. Also $%
\{\phi_{\lambda,c}\}$ converges pointwise to the constant function one as $%
\lambda$ tend to $0^+$ with fixed $c$.

Suppose we use this family of functions as a multiplier to determine a mean
for a variable obeying the Cauchy distribution. Let $m(\lambda,c)$ be
defined by:

\begin{equation*}
m(\lambda, c) = \int_{-\infty}^{\infty} \phi_{\lambda,c}(x)\frac{x}{\pi (1 +
x^2)} \, dx = \int_{-\infty}^{0} \frac{c \lambda e^{\lambda x} x^2}{ 1 + x^2}
\, dx = c m(\lambda,1). 
\end{equation*}

\noindent Now

\begin{eqnarray*}
1 = e^{\lambda x} |_{-\infty}^{0} & = & \int_{-\infty}^{0} \lambda
e^{\lambda x} \, dx \\
\geq \int_{-\infty}^{0} \frac{ \lambda e^{\lambda x} x^2}{ 1 + x^2} \, dx &
= & m(\lambda,1) = \int_{0}^{\infty}\frac{ \lambda e^{-\lambda x} x^2}{1 +
x^2} \, dx \\
\geq \int_{K}^{\infty} \frac{\lambda e^{-\lambda x} x^2}{1 + x^2} \, dx &
\geq & \frac{K^2 e^{-\lambda K}}{1 + K^2}
\end{eqnarray*}

\noindent for any positive real number $K$. Thus

\begin{equation*}
1 \geq \limsup_{\lambda \rightarrow 0+} m(\lambda, 1) \geq \liminf_{\lambda
\rightarrow 0+} m(\lambda, 1) \geq \frac{K^2}{1 + K^2}. 
\end{equation*}

\noindent Letting $K$ tend to infinity, we find that $\lim_{\lambda
\rightarrow 0+} m(\lambda,1) = 1$.

Hence the multiplier-induced mean of the standard Cauchy distribution is:

\begin{eqnarray*}
E_{mult}(X) = \lim_{\lambda \rightarrow 0+} \int_{-\infty}^{\infty}
\phi_{\lambda,c}(x)\frac{x}{\pi (1 + x^2)} \, dx & = & \lim_{\lambda
\rightarrow 0+} m(\lambda, c) \\
= \lim_{\lambda \rightarrow 0+} c m(\lambda, 1) & = & c \lim_{\lambda
\rightarrow 0+} m(\lambda,1) = c.
\end{eqnarray*}

However, $c$ was arbitrary depending on the choice of the multiplier!

Although some may consider the Cauchy distribution anomalous, we remind the
reader that its legitimacy and importance stem in part from the fact that it
is the quotient of two independent standard normal random variables. It has
application in physics under the name of the Lorentz distribution. Indeed
long-tailed and counter-intuitive distributions are increasingly important
in recent times (see Gumble \cite{Gum} or Taleb \cite{Tal}) in financial
mathematics, the study of natural and man-made disasters, and computer
network analysis. Extending the notion of mean to such distributions, and
investigating the limits of the notion of mean in such settings, are among
the ways of moving beyond the normal regime.

\section{State Space Theory}

State Space Theory or System Theory is widely used to provide a mathematical
description of physical systems including those of classical mechanics as
well as other systems such as biological and social systems. The state of
the system at any time is taken to be an element of a set $S$ called state
space. The evolution of the state is given by a function $T_{t}:S
\rightarrow S$ taking the state $s$ at time $0$ to the state $T_{t}(s)$ at
time t. \ A (real-valued) \emph{observable} is any function $f: S
\rightarrow R$ which assigns to each state $s$ a number $f(s)$ (see Mackey%
\cite{Mack}). All observables may be determined from the state, and indeed
the state can be viewed as a maximal independent set of observables that
characterize the system at a given time. The dynamic evolution of the state
is deterministic and time may be taken to be either discrete or continuous.
Evolution of an observable $f$ can be expressed by $t \mapsto f \circ
T_{t}(s)$, i.e., the value of the observable at time t is obtained by
applying the observable function to the state at time t.

A familiar example of the state space approach is Hamiltonian mechanics. The
state space in this case is phase space, and a state is a 2n-tuple \newline
$(q,p)=(q_{1},...,q_{n},p_{1},...,p_{n})$ consisting of position coordinates 
$q_{i}$ and momentum coordinates $p_{i}$. The evolution is $T_t(q,p) =
(q(t),p(t))$, where the latter is the solution to Hamilton's equations with
initial data $(q(0),p(0)) = (q,p)$:

\begin{eqnarray*}
\frac{dq_{i}}{dt} &=&\frac{\partial H}{\partial p_{i}} \\
\frac{dp_{i}}{dt} &=&-\frac{\partial H}{\partial q_{i}}
\end{eqnarray*}

\noindent for $i=1,2,...,n$. Here $H(q,p)$ is the Hamiltonian function of
the system, which is assumed to be a continuously differentiable function on
state space representing the total energy of the system. The function $H$ is
an example of an observable, as are the position and momentum coordinates,
angular momenta $q_{i}p_{j}-q_{j}p_{i}$, et cetera. A differentiable
observable $f$ evolves according to the equation:

\begin{equation*}
\frac{df}{dt} = \sum_{i = 1}^n (\frac{\partial f}{\partial q_i} \frac{%
\partial H}{\partial p_i} - \frac{\partial f}{\partial p_i} \frac{\partial H%
}{\partial q_i}), 
\end{equation*}

\noindent the right-hand side being the definition of the Poisson bracket $%
[f, H]$, under which operation $C^{\infty}$ observables form a Lie algebra.

Given a deterministic state space it is natural to pass to a statistical
setting as follows. We replace the old states $s$ by new states that are
(Borel) probability measures $P$ on the state space $S$. The old observables 
$f$ on the original state space are replaced by new observables that are the
means of the old observables with respect to the probability measure $P$.
Thus for any original state space observable $f$, the map

\begin{equation*}
P\mapsto E_{P}(f) 
\end{equation*}

\noindent defines an observable on the set of probability measures. If $f$
is bounded, this observable is defined for all probability measures. If not,
it is defined for those measures with respect to which $f$ is integrable.

If the only observables allowed were obtained in this manner, this would
appear to be a severe limitation. However, the variance of f and all moments
of f can themselves be regarded as means of original observables. Indeed,
even the probability distribution for $f$ can be regarded as a mean. This is
because on a Borel subset $A$ of the reals, the probability that $f$ takes a
value in $A$ is given by $E_{P}(\chi _{A}\circ f)$, where $\chi _{A}$ is the
characteristic function of the set $A$.

The evolution of the probabilistic state can be induced by an underlying
deterministic evolution. The probability measure at time $t$, $P_{t}$, is
given by $P_{t}(A)=P(T_{-t}(A))$ where $A$ is any (Borel) subset of $S$.
This permits us to talk about the evolution of observables since the mapping 
$t \mapsto E_{P_{t}}(f)$ describes such an evolution. In the Hamiltonian
formalism phase space has a natural 2n-dimensional Lebesgue measure $\lambda 
$ called Liouville measure with infinitesimal volume element $%
dq_{1}...dq_{n}dp_{1}...dp_{n}$, and $\lambda (T_{t}(A))=\lambda (A)$ for
all Borel subsets $A$ of $S$ and all times t. Dynamics in phase space can be
thought of as a fluid flow that permits change of shape but no change in
volume. The probability state $P$ can often be taken to be the integral of a
probability density function $\rho (q,p)$ with respect to $\lambda $. At the
other extreme $P$ can be taken to be a delta function $\delta
(q-q_{0})\delta (p-p_{0})$, which reduces to the deterministic theory with
state $s=(q_{0},p_{0})$. The probabilistic setting also permits us to
abandon the deterministic evolution $\{T_{t}\}$ and work with a stochastic
evolution exclusively, e.g., one of Markov type.

A use of means in state space theory that we have not touched on here
relates to ergodic theory, in which time averages of observables over
trajectories are compared with averages over state space regions using a
suitable normalized volume measure.

The essential point is that means provide the transition from classical
observables for deterministic systems to statistical observables for
stochastic systems.

\section{Entropy}

A subtlety occurs in statistical mechanics that is not present in ordinary
probability theory. An observable is commonly defined as a real-valued
function of the state, and in statistical mechanics the state is a
probability measure $P$ on state space. Thus any real-valued function of $P$
can be taken to be an observable, e. g., \ $P \mapsto P(B)$ is an observable
where $B$ is any fixed Borel set in the state space $S$. This observable is
an expected value since $P(B) = E_{P}(\chi _{B})$. However, not all
observables arise as expected values of original observables. The most
familiar example of such an observable is the entropy function, which can be
interpreted as an expected value (mean) but is not a conventional mean.

To avoid certain difficulties associated with the continuous case we will
confine our attention to the case where the underlying state space is a
finite set. Let $S$ be a finite state consisting of n states. A classical
observable is a function $f: S \rightarrow (\infty, \infty)$. A discrete
classical evolution might be a function $T: S \rightarrow S$ such that if $i$
is the state at a given time then $T(i)$ is the state one time unit later.
(A continuum of times presents a problem for deterministic evolution in a
finite state space, although that problem does not arise in the
probabilistic setting.)

When we pass to a statistical notion of state, we arrive at a probability
vector $p=(p_{1},...,p_{n})$ where $p_{i}$ is the probability that the
system is in the deterministic state $i$. We can now form expected values of
classical observables $f$, i.e.,

\begin{equation*}
E(f) = \sum_{i} p_{i}f(i) 
\end{equation*}

\noindent as noted before. We can also form such expressions as the entropy:

\begin{equation*}
H(p) = \sum_{i}p_{i}\log (\frac{1}{p_{i}}). 
\end{equation*}

\noindent Superficially the entropy appears to be another mean value, the
mean value of the ``uncertainty" $\log (1/p_{i})$, also called the
``surprise value" (The log here is usually taken with base $2$.) Thus the
entropy of a probability state is the mean uncertainty of the state. This is
not the mean of a classical observable since the function $i \mapsto \log
(1/p_{i})$ is not a classical observable. Classical observables should exist
and be measurable prior to assignment of probabilities, but it makes no
sense to consider the uncertainty function until probabilities have been
introduced. The dependency of the uncertainty function on $i$ is not
intrinsic and is only determined through the postulated probability state $%
p_i$.

It happens that entropy has another relationship to means of considerable
importance, namely through the \emph{Maximum Entropy Principle (MEP)}, also
known as \emph{Jaynes' Principle}. In the absence of an evolutionary law $T$
and an initial assignment, we are faced with the problem of determining the
probability state $p$, i. e., an assignment of probabilities to the
deterministic states $i$. The MEP \cite[p. 370]{Jaynes} asserts that:

\begin{quotation}
The probability state $p$ maximizing entropy subject to the given values $%
\alpha _{1},...\alpha _{k}$ for the means of known classical observables $%
g_{1},...,g_{k}$ provides predictions ``most strongly indicated by our
present information."
\end{quotation}

Using the calculus of variations, we can in general determine a unique
distribution among those that satisfy the constraints

\begin{equation*}
\sum_{i}p_{i}g_{j}(i)=\alpha _{j} 
\end{equation*}

\noindent for $j=1,...,k$ and maximizing $H(p)$, namely, the one with the
probability assignment

\begin{equation*}
p_{i}=C\exp \{-(\sum_{j}\beta _{j}g_{j}(i))\} 
\end{equation*}

\noindent for $i=1,...,n$ where $\beta_1, \dots, \beta_k$ are constants
determined from the $\alpha_i\mbox{'s}$, and $C$ is a positive normalizing
constant chosen so that the sum of the $p_{i}$'s is $1$.

The interpretation of this result takes two forms (at least). Suppose the
states are those of an individual particle in a gas of $N$ particles. Then
the quantities $\alpha _{1},...\alpha _{k}$ represent measured values of the
total value of $g_{1},...,g_{k}$\ over the entire gas divided by $N$. The
probabilities $p_{i}$, derived from the MEP, are the probabilities that a
particle picked at random from among the $N$ particles is in the i-th state.
They may also be regarded as the fraction of particles that are in the i-th
state. We may not care about individual particles but we do care about these
fractions, which can be taken to define the macroscopic state of the gas
(volume, pressure, temperature, and the like). This is the \emph{ensemble}
viewpoint of Gibbs. Yet another perspective is to regard what we usually
observe as a small perturbation about values induced by means.



\section{Quantum Issues}

The mean plays a pivotal role in quantum theory, even if this role has not
been examined closely in most treatments of quantum theory. In quantum
mechanics the state of a physical system is described by a wave function $%
\psi$ that is an element of a Hilbert space $\mathcal{H}$. (Strictly
speaking $\psi$ is not a function but an equivalence class of functions, and
in addition each state is associated with a ray in Hilbert space.) Each
physical observable that takes on real-number values (e.g., a position
coordinate, a momentum coordinate, the energy, a spin component) is
associate with a self-adjoint operator $A$ in $\mathcal{H}$. For simplicity
each observable is denoted by the same symbol ``A" as the associated
operator. Any self-adjoint operator A has in turn an associated
projection-valued measure $P_A$ (see, for example, \cite{Mack}) that assigns
to each Borel set $S$ in $\mathcal{R}$ an orthogonal projection $P_A(S)$ in
the Hilbert space:

\begin{equation*}
S \longmapsto P_A(S) 
\end{equation*}

\noindent in such a way that $A$ is an integral combination of these
orthogonal projections, represented symbolically by:

\begin{equation*}
A = \int_{\mathcal{R}} x \, P_A (dx), 
\end{equation*}

\noindent or by:

\begin{equation*}
A (\psi) = \int_{\mathcal{R}} x \, P_A (dx)(\psi), 
\end{equation*}

\noindent where $\psi$ is in the domain of $A$. If a measurement is made,
the probability that the value of A is in the set $S$ when the system state
is $\psi$ defined to be:

\begin{equation*}
\left\langle P_{A}(S)(\psi ),\psi \right\rangle =||P_{A}(S)(\psi )||^{2} 
\end{equation*}

\noindent where $<\,,\,>$ is the inner product on $\mathcal{H}$, linear in
the first variable and conjugate-linear in the second variable, and $||\,||$
is the norm on $\mathcal{H}$.

Quantum Mechanics is thus a statistical theory based on a family of
probability measures defined by:

\begin{equation*}
S \longmapsto ||P_A(S)(\psi)||^2. 
\end{equation*}

\noindent These are the Borel probability measures associated with
observables $A$ when the system state is $\psi$. One consequence of this is
that the set of possible values of $A$ is the spectrum of the operator $A$,
and another is that the mean of $A$, when the state is $\psi$, is given by:

\begin{equation*}
\left\langle A(\psi ),\psi \right\rangle =\int_{\mathcal{R}}x\,\left\langle
P_{A}(dx)(\psi ),\psi \right\rangle =\int_{\mathcal{R}}x||P_{A}(dx)(\psi
)||^{2}. 
\end{equation*}

\noindent In particular the quantity $\left\langle P_{A}(S)(\psi ),\psi
\right\rangle =||P_{A}(S)(\psi )||^{2}$ can be interpreted as the mean of
the observable $P_{A}(S)$ when the state is $\psi $. The observable $%
P_{A}(S) $ is an orthogonal projection, taking the value $1$ when the value
of $A$ is in is $S$ and the value $0$ when the value of $A$ is not in $S$.
Thus $||P_{A}(S)(\psi )||^{2}$ also represents the probability that $A$ is
in $S$ when the state is $\psi $. This is a reminder that all probabilities
are means.

The mean, $\left\langle A(\psi ),\psi \right\rangle $, is the integral over
the real line of the real variable $x$ with respect to the Borel probability
measure $||P_{A}(\,)(\psi )||^{2}$. Thus the mean exists, it appears, if and
only if $x$ is integrable with respect to this measure, thus if and only if $%
\psi $ is in the domain of $A$. Self-adjoint operators have domains that are
dense in $\mathcal{H}$ but many of the most prominent ones (e.g., those
associated with position and momentum and often energy) do not have domain
equal to $\mathcal{H}$. Hence there will be states for which the means of
some observable may not be well-defined. Whether these states are realizable
in practice is uncertain, but there is no good theoretical reason why they
should be ignored. (Our discussion focuses on mathematical definition and
characterization. The spectrum of a self-adjoint operator is identified with
the possible values of a measured quantity. If the spectrum is discrete, a
measurement may be able to distinguish one value from another; if the
spectrum is continuous, measurement will only be able to determine an
interval that contains the value, not the exact value. Repeated measurements
when the system is in the same state thus only arrive at a rough
approximation of the distribution and a rough estimate of the mean for a
state.)

A curiosity in quantum mechanics, not ordinarily seen in other applications
of probability, is the following. Suppose $\mu =\left\langle A(\psi ),\psi
\right\rangle $ is the mean of some observable $A $ when the state is $\psi $%
. Then the variance of the observable in this state is naturally given by:

\begin{equation*}
\int_{\mathcal{R}}(x-\mu )^{2}\,\left\langle P_{A}(dx)(\psi ),\psi
\right\rangle \, = \,\left\langle (A-\mu I)^{2}(\psi ),\psi \right\rangle
=||(A-\mu I)(\psi )||^{2}. 
\end{equation*}

\noindent So the variance exists if and only if $\psi$ is in the domain of $A
$. The condition for the mean to exist is the same as the condition for the
variance to exist. In quantum mechanics we are led to think that the only
distributions for which the mean is finite are ones in which the variance is
also finite. However, a closer look at this situation reveals some
discrepancies.

The chief discrepancy is the following. Suppose that the original observable 
$A$ can be written in the form

\begin{equation*}
A = C - D 
\end{equation*}

\noindent where $C$ and $D$ are non-negative self-adjoint operators.
Non-negative self-adjoint operators can be written as squares of
self-adjoint operators, so that $C = E^2$ and $D = F^2$ with $E$ and $F$
self-adjoint. Then

\begin{eqnarray*}
\left\langle A(\psi ),\psi \right\rangle =\left\langle (C-D)(\psi ),\psi
\right\rangle &=&\left\langle (E^{2}-F^{2})(\psi ),\psi \right\rangle \\
=\left\langle E^{2}(\psi ),\psi \right\rangle -\left\langle F^{2}(\psi
),\psi \right\rangle &=&||E(\psi )||^{2}-||F(\psi )||^{2}.
\end{eqnarray*}

\noindent Thus, the mean of $A$ exists if and only if $\psi$ is in the
intersection of the domain of $E$ and the domain of $F$. It is easy,
incidentally, to construct examples of elements of $\mathcal{H}$ that are in
the domain of a self-adjoint operator $E$ but are not in the domain of its
square $E^2$. In addition, as it happens, it is possible to offer explicit
candidates for the operators $E$ and $F$ given $A$. Set

\begin{equation*}
E = \int_{(0, \infty)} \sqrt{x} \, P_A(dx) \mbox{ and } F = \int_{(-\infty,
0)} \sqrt{-x} \, P_A(dx). 
\end{equation*}

\noindent We conclude that the mean of $A$ exists when $\psi$ is in $%
\mbox{dom }E \cap \mbox{dom } F$ and the variance of $A$ exists when $\psi$
is in $\mbox{dom } A$. If $\psi$ is not in $\mbox{dom } E$ but is in $%
\mbox{dom }F$, then it is reasonable to say that the mean of $A$ is $\infty$%
. Likewise if $\psi$ is in $\mbox{dom }E$ but not in $\mbox{dom }F$, the
mean of $A$ is $-\infty$. If $\psi$ is in neither $\mbox{dom }E$ nor $%
\mbox{dom }F$, then the ordinary mean does not exist. In the spirit of our
discussion of $L(c)$ earlier, it is possible to truncate the integrals for $E
$ and $F$ in the last display, replacing $\infty$ by $M$ and $-\infty$ by $-M
$ and investigate the existence of an appropriate combined limit as $M$ tend
to infinity.

Similar considerations can be applied when the ``pure" state $\psi $ is
replaced by a density matrix representing a statistical ensemble of pure
states, or in rigged Hilbert spaces where the existence of states varies
according to the properties of the observable, or to cases arising by the
use of positive-operator-valued measures generalizing the projection-valued
measures treated above.


\section{Conclusion}

The mean, as we have seen, is ubiquitous in scientific explanation. Not only
does it provide a summary of sample data and, when it exists, of data from
the entire population, but it establishes a connection between samples and
the whole population. Furthermore, it facilitates generalization of
deterministic observables that are functions of the deterministic state to
probabilistic observables that are functions of the probabilistic state. The
Maximum Entropy Principle then makes use of constrained means to identify
macroscopic distribution of physical importance parametrized by these means.
While quantum mechanics abandons determinism, it retains the notion of mean
to summarize the possible results of experiments and the measurement of
quantum observables. Although not all observables have finite means, weak
and doubly weak means and the alternatives identified in Theorem 3.1 provide
an enumeration of possible behaviors of variables and associated probability
distributions, and give further insight into potentialities associated with
large data sets.

\end{document}